\documentclass[12pt]{amsart}
\usepackage{amsmath,amscd,amssymb}
\usepackage[arrow,matrix]{xy}

\theoremstyle{plain}
\newtheorem{theorem}[subsection]{Theorem}
\newtheorem{lemma}[subsection]{Lemma}

\theoremstyle{definition}
\newtheorem{remark}[subsection]{Remark}

\newtheorem{example}[subsection]{Example}

\numberwithin{equation}{section}
\newcommand{\A}{{\mathcal A}}
\newcommand{\I}{{\mathcal I}}
\newcommand{\B}{{\mathcal B}}
\newcommand{\CC}{{\mathcal C}}
\newcommand{\QQ}{{\mathcal Q}}

\newcommand{\Z}{\mathbb{Z}}
\newcommand{\Q}{\mathbb{Q}}
\newcommand{\R}{\mathbb{R}}
\newcommand{\C}{\mathbb{C}}

\newcommand{\N}{\mathbb{N}}
\newcommand{\PP}{\mathbb{P}}

\DeclareMathOperator{\diag}{diag}
\DeclareMathOperator{\sgn}{sgn}

\newcommand{\eqv}{{\Longleftrightarrow}}

\begin{document}

\title[Manifolds from complete intersections]{Closed manifolds 
coming from\\
Artinian complete intersections}

\author[\c{S}tefan Papadima]{\c{S}tefan Papadima$^*$}
\address{Institute of Mathematics "Simion Stoilow",
P.O.Box 1-764,
RO-014700 Bucharest, Romania}
\email
{Stefan.Papadima@imar.ro}

\author[Lauren\c{t}iu P\u{a}unescu]
{Lauren\c{t}iu P\u{a}unescu$^\dagger$}
\address{School of Mathematics and Statistics,
University of Sydney, Sydney,
New South Wales 2006, Australia}
\email
{laurent@maths.usyd.edu.au}

\thanks{$^*$Partially supported by grant U4249 Sesqui R\&D/2003
of the University of Sydney.} 

\thanks{$^\dagger$Partially supported by grant U4249 Sesqui R\&D/2003
of the University of Sydney.}



\begin{abstract}
We reformulate the integrality
property of the Poincar\'{e} inner product in the middle dimension,
for an arbitrary Poincar\'{e} $\Q$-algebra, in classical terms
(discriminant and local invariants). When the algebra is $1$-connected,
we show that this property is the only obstruction to realizing it
by a closed manifold, up to dimension $11$. 
We reinterpret a result of Eisenbud and Levine on finite map germs, 
relating the degree of the map germ to the signature of the 
associated local ring, to
answer a question of Halperin on artinian weighted 
complete intersections.
We analyse the homogeneous artinian
complete intersections over $\Q$ realized by closed manifolds of
dimensions $4$ and $8$, and their signatures. 
\end{abstract}

\maketitle

\section{Introduction}
\label{sec:intro}

\subsection{Artinian complete intersection.}
\label{subsec=ci}
Let $\A$ be a {\it weighted artinian complete intersection}
($WACI$), that is, a commutative graded $\Q$-algebra of the form
\begin{equation}
\label{eq:WACI}
\A = \Q[x_1, \dots, x_n]/ \I \, ,
\end{equation}
where the variables $x_i$ have positive even weights, 
$w_i := \mid x_i \mid$, 
and the ideal $\I$ is generated by a regular sequence,
\begin{equation}
\label{eq:reg}
\I = (f_1, \dots, f_n) \, ,
\end{equation}
of weighted-homogeneous polynomials, $f_i$.

One knows~\cite[Theorem 3 and p.198]{H} that $\A^*$ is a
$1$-connected {\it rational Poincar\' e duality algebra}
($\Q$--$PDA$), with Poincar\' e polynomial
\begin{equation}
\label{eq=fdim}
\A^* (t)= \prod_{i=1}^n \frac{1-t^{\mid f_i\mid}}{1- t^{\mid x_i \mid}}
\end{equation}
(and, consequently, with even formal dimension,
$m =\sum_{i=1}^n (\mid f_i\mid -\mid x_i\mid)$).

\subsection{The integrality obstruction}
\label{subsec:int}

Let $\A^*$ be an arbitrary $1$-connected Poincar\' e
duality $\Q$-algebra, with formal dimension $m$. The
{\it smoothing problem} we are going to look at is the
following:
\begin{quote}
is $\A^*$ isomorphic to a graded algebra of the form
$H^*(M^m, \Q)$, where $M$ is a $1$-connected closed
smooth $m$-manifold?
\end{quote}
We shall say that $\A$ is {\em smoothable} if the answer is yes.

By $\Q$-surgery (\cite{S}, \cite{Ba}), we know that $\A$ is
smoothable, for $m\ne 4k$.
Assume now that $m=4k$, and pick an {\it orientation},
$\omega \in \A^{4k}\setminus \{ 0 \}$. This gives rise
(via Poincar\' e duality) to a symmetric inner product space
over $\Q$, denoted by $(\A^{2k}, \cdot_{\omega}) \in W(\Q)$.
(Here and in the sequel, $W(R)$ denotes the Witt group of the
ring $R$; see~\cite{MH}.)
If $\A$ is smoothable, then clearly
\begin{equation}
\label{eq:z}
(\A^{2k}, \cdot_{\omega}) \in W(\Z), \quad
\text{for some orientation} \quad \omega \, .
\end{equation}

It turns out that the {\it integrality obstruction} from \eqref{eq:z},
for a fixed orientation $\omega$, 
is equivalent to the fact that the quadratic form 
on $\A^{2k}$ associated to $\cdot_{\omega}$ 
is a sum of signed squares, over $\Q$; 
see~\cite[Corollary IV.2.6]{MH}.

When the signature is zero, the integrality condition is 
equivalent to $(\A^{2k}, \cdot)$ being split; 
see \cite[I.6--7]{MH}. In this case,
\eqref{eq:z} is the only obstruction to smoothing; see
\cite{S} and \cite{Ba}, and also 
\cite[Proposition 3.4]{Pp}. In the 
non-zero signature case, additional obstructions 
may appear, see e.g.~\cite[\S 4.5]{Pp} for some simple
examples, based on \cite{Ba} and~\cite{BLLV}.

\subsection{Main results}
\label{subsec=res}

The smoothing problem described above may be solved by using
fundamental $\Q$--surgery results due to D.~Sullivan; see
\cite{S}, \cite{Ba}. This opens the way for constructing 
closed manifolds with interesting geometric properties, starting
from $\Q$--$PDA$'s (see~\cite{Pp} for applications to geodesics).

The difficulties of the smoothing problem stem from the fact that 
the obstructions involve, besides \eqref{eq:z}, delicate conditions 
on the signature. If $\A$ is an arbitrary $WACI$, we point out 
a topological interpretation of the signature, in terms of the
defining relations of $\A$, thus answering a question of
S.~Halperin from~\cite{H}. We do this in Section~\ref{sec:sign},
by using a basic result on finite map germs, due to D.~Eisenbud
and H.~Levine~\cite{EL}.

In Section~\ref{sec:low}, we focus on the integrality obstruction.
In Theorem~\ref{thm:eight}, we show that \eqref{eq:z} is the only
obstruction to smoothing, for arbitrary $\Q$--$PDA$'s of formal
dimension $4$ or $8$. This is no longer true in dimension $12$; 
see Remark~\ref{rem:best}. In dimension $8$, our proof requires a 
classical result on sums of four squares.

As far as condition~\eqref{eq:z} is concerned, it may be handled,
for a fixed orientation, by using discriminants and local
invariants of nondegenerate quadratic forms over $\Q$; see~\cite{Se}.
We give a similar interpretation for~\eqref{eq:z}, where $\A$ is
an arbitrary $\Q$--$PDA$, in Theorem~\ref{thm=test}, by analysing
changes of orientation. For odd rank, the answer depends only on 
local invariants. For even rank, both local invariants and
discriminant are involved, in general; see Remark~\ref{rk=w}.

The results from Section~\ref{sec:low} are applied in
Section~\ref{sec:hom}. Here, we construct $8$--manifolds with
interesting properties, starting from $WACI$'s which are
homogeneous (that is, with $w_i=2$, for all $i$). The integrality
test from Theorem~\ref{thm=test} is illustrated on two families
of examples: one with odd rank (see Example~\ref{ex=c1}), and the
other with even rank (see Example~\ref{ex:par}). The even rank family 
has the remarkable property that the corresponding test, described
in Theorem~\ref{thm=test}~\eqref{eq=t2}, collapses to a single,
simple, discriminant obstruction.

\section{Signature and degree}
\label{sec:sign}

The signature plays an important role in the smoothing problem
described in \S~\ref{subsec:int}, via the Hirzebruch formula; 
see \cite{M}, \cite{S}, \cite{Ba}.

Let $\A= \Q[x_1, \dots, x_n]/(f_1, \dots, f_n)$ be 
an arbitrary $WACI$, as defined in  \S \ref{subsec=ci}.
Among other things, S. Halperin showed in \cite[Theorem 3]{H}
that $\A^*$ is a Poincar\' e duality algebra ($1$-connected and
commutative), giving thus rise to $(\A, \cdot_{\omega})\in W(\Q)$,
for any choice of orientation, $\omega\in \A^m\setminus \{ 0 \}$.
At the end of section \S $9$ from \cite{H}, he raised 
the following question: is there an explicit way of computing
the signature, $\sigma (\A, \cdot_{\omega})$, in terms of the 
defining relations of $\A$?

Let us consider the associated $C^{\infty}$ map germ,
\begin{equation}
\label{eq=germ}
f: (\R^n, 0) \rightarrow (\R^n, 0) \; ,
\end{equation}
having as components the defining polynomial relations of $\A$.
Since $f^{-1}(0) = \{ 0 \}$, the degree at $0$ of
$f$, $\deg\, (f)$, is defined, and may be computed in terms of
regular values of $f$.

Using Theorem 1.2 from Eisenbud and Levine~\cite{EL}
(see also~\cite[p.103--104]{A} and~\cite{K}), we may offer 
the following answer to the above signature problem.

\begin{theorem}
\label{thm:sign}

Let $\A$ be an arbitrary $WACI$. 
Let $f: (\R^n, 0) \rightarrow (\R^n, 0)$ be the 
$C^{\infty}$ map germ associated to a system of 
defining relations for $\A$. Then:
\[
\sigma(\A, \cdot_{\omega}) = \deg\, (f)\; ,
\]
for a good choice of orientation, $\omega$.
\end{theorem}

\begin{proof}
In \cite {EL}, the authors give a concrete way to algebraically compute 
the degree of a {\em finite} map.
Namely let us consider $\QQ(f)= \CC_0^{\infty}(\R^n)/(f_1, \dots, f_n)$, 
where  $\CC_0^{\infty}(\R^n)$ is the ring of germs 
at $0$ of smooth real-valued functions. A map 
$f=(f_1, \dots, f_n):(\R^n, 0) \rightarrow (\R^n, 0) $ is 
called finite exactly when the corresponding $\QQ(f)$ is a finite 
dimensional real vector space. Consequently $0$ is isolated in $f^{-1}(0)$ 
and this  allows one to consider the topological degree of $f$.
In this set-up  (\cite[Theorem 1.2]{EL}) one can compute the degree of $f$ 
as the signature of a symmetric bilinear form on $\QQ(f)$.
More precisely, if we define $\omega$ to be the class of the Jacobian 
$J =\det\,(\partial f_i/\partial x_j)$ in $\QQ (f)$, 
it follows from~\cite[Proposition 4.4(ii)]{EL} that $\omega\ne 0$. Pick
any linear functional, $\varphi : \QQ(f) \to \R$, 
such that $\varphi (\omega) > 0$. A symmetric bilinear form on $\QQ(f)$ 
whose signature gives the degree is given 
by the following formula:
\begin{equation}
\label{eq=inn}
<p,q>_{\varphi}= \varphi (pq)\, , \quad \text{for} \quad p,q \in \QQ(f)\, . 
\end{equation}

This result may be applied to the $C^{\infty}$ map germ \eqref{eq=germ},
in the following way. Firstly, note that the natural algebra map
\begin{equation}
\label{eq=iso1}
\Phi \colon \A \otimes \R = \R [x_1,\dots ,x_n]/(f_1, \dots ,f_n)
\rightarrow \R[[x_1,\dots,x_n]]/(f_1, \dots, f_n)
\end{equation}
is an isomorphism. Indeed, the surjectivity of $\Phi$ easily follows from
the fact that $(x_1, \dots ,x_n)^{m+1} \subset (f_1, \dots ,f_n)$
(guaranteed by $\A^{>m} =0$), while injectivity may be checked 
by a straightforward weight argument.

We may now use the isomorphism \eqref{eq=iso1}
to show that the natural algebra map
\begin{equation}
\label{eq=iso2}
\Psi \colon  \A \otimes \R \rightarrow \QQ (f)
\end{equation}
is an isomorphism as well (and consequently $f$ is a finite map germ).
To this end, consider the algebra map
\[
\tau \colon \CC_0^{\infty}(\R^n)/(f_1, \dots, f_n)
\rightarrow \R[[x_1,\dots,x_n]]/(f_1, \dots, f_n) \, ,
\]
induced by Taylor series expansion. Since $\tau \circ \Psi =\Phi$,
$\Psi$ must be injective. The surjectivity of $\Psi$ follows
from the well-known fact (see e.g. \cite[p.2]{m}) that
the ideal of $k$-flat functions (i.e., germs in $\CC_0^{\infty}(\R^n)$
whose derivatives vanish at $0$, up to order $k$) is contained in
$(x_1, \dots ,x_n)^{k+1}$, for any $k\ge 0$. 

We infer from \eqref{eq=fdim} and the isomorphism \eqref{eq=iso2}
that the class of $J$ modulo $\I$, $\omega$, belongs to 
$\A^m \setminus \{ 0 \}$. This is the desired orientation.
Indeed, we may use $\Psi$ to identify the algebras $\A \otimes\R$
and $\QQ (f)$, and then define $\varphi$ 
to be the linear projection from the graded vector space $\A^* \otimes\R$
to $\A^m \otimes\R \equiv \R$. Since $\varphi (\omega)=1$, 
the Eisenbud--Levine formula applies.

Obviously, the bilinear form \eqref{eq=inn} on $\QQ(f)$ is identified
with the Poincar\' e inner product on $\A \otimes \R$,
$\cdot_{\omega}$. Our signature formula follows. 
\end{proof}

\section{Smoothing in small dimensions, and the integrality condition}
\label{sec:low}

\subsection{Small dimensions}
\label{subsec=low}

We begin by showing how the general smoothing problem
becomes simpler, in dimensions up to $11$. 

\begin{theorem}
\label{thm:eight}

Let $\A$ be a $1$-connected Poincar\' e duality $\Q$-algebra
with formal dimension $m$ (not necessarily a $WACI$). Assume
that $m\le 11$. Then: $\A$ is smoothable (in the sense 
explained in {\em \S \ref{subsec:int}}) if and only if
there is $\omega \in \A^{4k}\setminus \{ 0 \}$ such that 
$(\A^{2k}, \cdot_{\omega})\in W(\Z)$, when $m=4k$, and $\A$ is
always smoothable, otherwise.
\end{theorem}

\begin{proof}
For $m \not\equiv 0$ (mod $4$) arbitrary, smoothability
follows from \cite{S} and \cite{Ba}. Assume then that
$m=4k$, with $k=1$ or $2$.

We have to show that $\A$ is smoothable, as soon as 
property \eqref{eq:z} from \S \ref{subsec:int} holds.
Pick an orientation $\omega$ such that 
$\sigma := \sigma (\A^{2k}, \cdot_{\omega}) \ge 0$.
We know that the Poincar\' e quadratic form on
$\A^{2k}$ is a sum of $t$ squares, $t\ge \sigma$,
minus a sum of $s$ squares. 

If $k=1$, then plainly 
$\A^* =H^*((\#_t \C \PP^2)\# (\#_s \overline{\C \PP^2}), \Q)$,
for $t+s>0$, and $\A^* = H^*(S^4, \Q)$, for $t+s=0$.

If $k=2$ and $\sigma =0$, smoothability is guaranteed by 
\eqref{eq:z}; see \S \ref{subsec:int}. 

Assume then that $k=2$ and $\sigma >0$. We claim that if
the system
\begin{equation}
\label{eq:sq}
\begin{cases}
a+b &= \sigma \\
25a+18b & = \sum^t_{i=1} \alpha^2_i
\end{cases}
\end{equation}
has integer solutions, then $\A$ is smoothable.

Indeed, we may take the following algebraic Pontrjagin classes:
$q_2 = (10a+9b) \omega$, and
$q_1 = \sum^t_{i=1} \alpha_i x_i$, where
$\{ x_i \}$ is the canonical basis of the positive
definite part of $\A^4$. Set 
$N^8 = a\cdot \C \PP^4 + b\cdot \C \PP^2 \times \C \PP^2$.
Using the second equation from \eqref{eq:sq},
one may easily check that 
$\A$ and $N$ have the same Pontrjagin numbers;
see \cite{M}. This implies, via the first equation from
\eqref{eq:sq}, that the Hirzebruch signature formula
holds for $\A$, and we are done (\cite{S}, \cite{Ba}).

We come back to the system \eqref{eq:sq}. If $t\ge 4$,
the theorem of Bachet de M\' eziriac--Lagrange~\cite[II.8]{MH}
guarantees integer solutions.

In the remaining cases, $\sigma$ must be $1$, $2$ or $3$,
and then \eqref{eq:sq} may be solved as follows:
\[
\left\{
\begin{array}{llclcl}
a=1 \, , & b=0 & \text{and} & 25=5^2\, ,& \text{for}& \sigma =1\, ;\\
a=0 \, , & b=2 & \text{and} & 36=6^2\, ,& \text{for}& \sigma =2\, ;\\
a=1 \, , & b=2 & \text{and} & 61=5^2+6^2 \, ,& \text{for}& \sigma =3\, .
\end{array}
\right.
\]
This completes our proof.
\end{proof}

\begin{remark}
\label{rem:best}
The range $m\le 11$ is the best one for which 
the integrality condition alone guarantees smoothability.
Indeed, $\A = \Q[x]/(x^3)$, with $\mid x \mid =6$ has the
integrality property, without being smoothable. 
See~\cite[\S 4.5]{Pp}.
\end{remark}

\subsection{An integrality test}
\label{subsec=intor}

In applications, we will need to check the integrality condition
\eqref{eq:z}. This is clearly related to the theory of
nondegenerate quadratic forms over $\Q$. We thus start by 
reviewing some relevant facts from~\cite{Se}.

To begin with, assume that $\A$ is an arbitrary $\Q$--$PDA$,
with formal dimension $4k$. Set $r:= \dim_{\Q} \A^{2k}$,
and choose a $\Q$--basis of $\A^{2k}$. Pick any orientation,
$\omega \in \A^{4k}\setminus \{ 0 \}$, and denote by
$A_{\omega}$ the matrix of $\cdot_{\omega}$. Note that
$A_{\lambda \omega} =\lambda^{-1}\cdot A_{\omega}$, for any
$\lambda\in \Q^*$. The condition 
$(\A^{2k}, \cdot_{\omega}) \in W(\Z)$
translates to the fact that $A_{\omega}$ is equivalent
over $\Q$ (in the classical sense, see~\cite[IV.1]{Se}) with
a diagonal matrix of signs.

By a convenient choice of basis of $\A^{2k}$, we may suppose 
that $A_{\omega} =\diag (a_1,\dots, a_r)$. Then the discriminant
of $\cdot_{\omega}$ is equal to $a_1\cdots a_r$ (modulo $\Q^{*2}$).
For each prime number $p$, one also has a {\em local invariant}
at $p$, denoted by 
$$\varepsilon_p (A_{\omega}) := \prod_{1\le i<j \le r}
(a_i, a_j)_p \, \in \{ \pm 1 \} \, ,$$
where $(\cdot, \cdot)_p$ denotes the $p$--adic Hilbert symbol.
From the classification theory (\cite[IV.3]{Se}), we infer that
\begin{equation}
\label{eq=clas}
(\A^{2k}, \cdot_{\omega}) \in W(\Z)\eqv 
\mid a_1\cdots a_r \mid \in \Q^{*2}\quad
\text{and} \quad \varepsilon_p (A_{\omega})=1,
\forall p\equiv 1(2)\, .
\end{equation} 

Our next result translates in similar terms condition
\eqref{eq:z}, by taking into account changes of orientation.

\begin{theorem}
\label{thm=test}

Let $\A$ be an arbitrary $\Q$--$PDA$, of formal dimension $4k$. Then:
\begin{enumerate}
\item \label{eq=t1}
Assume $r\equiv 1(2)$. For any orientation $\omega$, there is
$\lambda \in \Q^*$ such that the discriminant of
$\cdot_{\lambda \omega}$ is $1$ (modulo $\Q^{*2}$). Supposing that
$\omega$ has the property that $a_1\cdots a_r \in \Q^{*2}$, 
\eqref{eq:z} is equivalent to
\[
\varepsilon_p (A_{\omega})=1,\,
\forall\, p\equiv 1(2)\, .
\]
\item \label{eq=t2}
Assume $r\equiv 0(2)$. Let $\omega$ be an arbitrary orientation.
Set $\epsilon := \sgn (a_1\cdots a_r)$. Two cases may occur:
\begin{itemize}
\item[$(\bullet)$]
$r\equiv 0(4)$ and $\epsilon =+1$, or 
$r\equiv 2(4)$ and $\epsilon =-1$. In this case, 
\eqref{eq:z} is equivalent to
\[
\mid a_1\cdots a_r \mid \in \Q^{*2}\quad
\text{and} \quad \varepsilon_p (A_{\omega})=1,
\forall p\equiv 1(2)\, .
\]
\item[$(\bullet \bullet)$]
$r\equiv 0(4)$ and $\epsilon =-1$, or 
$r\equiv 2(4)$ and $\epsilon =+1$. In this case, 
\eqref{eq:z} is equivalent to
\[
\mid a_1\cdots a_r \mid \in \Q^{*2}\quad
\text{and} \quad \varepsilon_p (A_{\omega})=1,
\forall p\equiv 1(4)\, .
\]
\end{itemize}
\end{enumerate}
\end{theorem}

\begin{proof}
Part \eqref{eq=t1}. The first assertion is easy: one may
take for instance $\lambda =a_1 \cdots a_r$. As for the second one,
it will follow from \eqref{eq=clas}, as soon as the following
claim is proved: if $(\A^{2k} , \cdot_{\lambda \omega})\in W(\Z)$,
where $\lambda >0$, then $(\A^{2k} , \cdot_{\omega})\in W(\Z)$.
To verify this claim, note that the first condition in \eqref{eq=clas}
implies that necessarily $\lambda \in \Q^{*2}$; this in turn
ensures that $\varepsilon_p (\lambda^{-1} A_{\omega})=
\varepsilon_p (A_{\omega})$, for all $p$ (since Hilbert symbols
are well-defined modulo squares, see~\cite[III.1]{Se}), and
we are done.

Part \eqref{eq=t2}. If $r$ is even, it readily follows 
from \eqref{eq=clas} that \eqref{eq:z} is equivalent to
the fact that there is $\lambda \in \Q^*$, $\lambda >0$, 
having the property that
\begin{equation}
\label{eq=sel}
\mid a_1\cdots a_r \mid \in \Q^{*2}\quad
\text{and} \quad \varepsilon_p (\lambda A_{\omega})=1,
\forall p\equiv 1(2)\, .
\end{equation}

It remains to compute the local invariants of
$\lambda \cdot A_{\omega}$, at odd primes. 
This may be done as follows. Firstly, one may use the
bilinearity of Hilbert symbols (\cite[III.1]{Se}),
together with the first property from \eqref{eq=sel},
to see that
\begin{equation}
\label{eq=loc}
\varepsilon_p (\lambda A_{\omega})=\varepsilon_p (A_{\omega})
\cdot (\lambda ,\lambda)_p^{\frac{r(r-1)}{2}}\cdot 
(\lambda ,\epsilon)_p \; .
\end{equation}

$(\bullet)$ In this case, elementary properties of Hilbert symbols
(\cite[III.1]{Se}) imply that \eqref{eq=loc} above reduces to
\begin{equation}
\label{eq=loc1}
\varepsilon_p (\lambda A_{\omega})=\varepsilon_p (A_{\omega})\, ,
\end{equation}
and we are done.

$(\bullet \bullet)$ Similarly, in this case \eqref{eq=loc} becomes
\begin{equation}
\label{eq=loc2}
\varepsilon_p (\lambda A_{\omega})=\varepsilon_p (A_{\omega})
\cdot (\lambda ,\lambda)_p \, .
\end{equation}
Note that $(\mu \nu ,\mu \nu)_p=(\mu ,\mu)_p (\nu ,\nu)_p$ and
$(2,2)_p=1$ (\cite[III.1]{Se}). It follows that we may assume in
\eqref{eq=sel} that $\lambda$ is a product of distinct odd
primes, $\lambda =q_1\cdots q_l$. Use~\cite[III.1]{Se} to compute
\begin{equation}
\label{eq=ll1}
(\lambda ,\lambda)_p= \prod_{i=1}^l (q_i ,q_i)_p=
\begin{cases}
1\; , & \text{for} \quad p\ne q_1,\dots ,q_l\, ;\\
(-1)^{\varepsilon (q_j)}\, , & \text{for} \quad p=q_j \, ,
\end{cases}
\end{equation}
where $\varepsilon(q)$ denotes the residue class modulo $2$
of $\frac{q-1}{2}$, as in~\cite{Se}. We infer from \eqref{eq=sel}, 
\eqref{eq=loc2} and \eqref{eq=ll1} that \eqref{eq:z} implies the
conditions from our statement.

Conversely, set
\begin{equation}
\label{eq=ll2}
\{ p=\, \text{odd}\, \mid \, \varepsilon_p (A_{\omega})=-1 \}=
\{ q_1, \dots, q_l \}\, .
\end{equation}
If all primes $q_j$ appearing in \eqref{eq=ll2} above are
equal to $3$ (modulo $4$), then we may take 
$\lambda =q_1 \cdots q_l$, and \eqref{eq:z} follows, again
from \eqref{eq=sel}, \eqref{eq=loc2}, and \eqref{eq=ll1}.
Our proof is complete.
\end{proof}

\begin{remark}
\label{rk=w}

Let $(V, \cdot)$ be a symmetric inner product space over $\Q$
(alias, a nondegenerate quadratic $\Q$--form). For any $k$, 
$(V, \cdot)$ may obviously be realized as $(\A^{2k}, \cdot_{\omega})$,
where the oriented $\Q$--$PDA$ $\A^*$ is $\Q \cdot 1$, in degree
$*=0$, $V$ in degree $*=2k$, $\Q \cdot \omega$ in degree $*=4k$, and
$0$ otherwise, with product given by $\cdot$.

Note first that the condition $(\A^{2k}, \cdot_{\lambda \omega})\in W(\Z)$
may depend on $\lambda \in \Q^*$. For odd $r$, examples are easy 
to construct, using the discriminant obstruction from \eqref{eq=clas}.
For $r=2$, for instance, a simple example is provided by the
quadratic form with matrix $A=\diag (5,5)$. Here, 
$(\A^{2k}, \cdot_{\omega})\in W(\Z)$, while
$(\A^{2k}, \cdot_{3 \omega})\notin W(\Z)$, even though the discriminant
condition from \eqref{eq=clas} is verified.

Note also that, in general, the local invariants show up in an 
essential way, in our integrality test from Theorem~\ref{thm=test}
Part~\eqref{eq=t2}. Indeed, consider the matrix
$A=\diag (1,1,1,2,5,10)$. The associated $\Q$--$PDA$ belongs to
case $(\bullet \bullet)$, and satisfies the discriminant condition
therefrom. On the other hand, $\varepsilon_5(A)=-1$, as readily seen.
\end{remark}

\section{Homogeneous complete intersections}
\label{sec:hom}

In this section, we want to apply Theorem~\ref{thm:eight}
to $WACI$'s, in the nontrivial cases, that is, when the formal
dimension is $4$ or $8$. We will restrict our attention to
{\it homogeneous} $WACI$'s, i.e., those with
$w_i =2$ and $\mid f_i \mid =2d_i \ge 4$, for
all $i$. Both conditions are very natural. The first one 
simply means that each $f_i$ is a homogeneous polynomial
of degree $d_i$. The restrictions $d_i \ge 2$ $(1\le i\le n)$
are imposed to avoid unnecessary redundancies, like
$\Q [x]/(x)=\Q$.

It is straightforward to check that
the formal dimension is equal to $4k$, with $k=1$ or $2$,
precisely in the cases listed below
(where $\underline{d}$ denotes $(d_1, \dots ,d_n)$, and
$r:= \dim_{\Q} \A^{2k}$); see \eqref{eq=fdim}.

Formal dimension $4$:

\[
\begin{array}{llr}
(I_4) & \underline{d} = (2,2)\, ; & r=2\, .\\
(II_4) & \underline{d} = (3)\, ; & r=1\, .
\end{array}
\]

Formal dimension $8$:

\[
\begin{array}{llr}
(I_8) & \underline{d} = (2,2,2,2)\, ;& r=6\, .\\
(II_8)  & \underline{d} =(2,2,3)\, ; & r=4\, .\\
(III_8)  & \underline{d} = (2,4)\, ; & r=2\, .\\
(IV_8) & \underline{d} = (3,3)\, ; & r=3\, .\\
(V_8) &  \underline{d} = (5)\, ; & r=1\, .
\end{array}
\]

Note that, in the general homogeneous $WACI$ case,
every degree vector, $\underline{d}$, may be 
realized by a smooth manifold. Indeed,
$H^*(\prod^n_{i=1} \C \PP^{d_i -1}, \Q)= 
\otimes^n_{i=1} \Q[x_i]/(x^{d_i}_i)$, with signature
$1$, when all $d_i$'s are odd, and $0$, otherwise. One may
ask whether more interesting signatures may also arise
from smooth manifolds. For instance, in case $(I_8)$ above,
the possible (non-negative) values of the signature are
$0,2,4,6$ (since $r=6$). Our last main result completely 
clarifies this question.

\begin{theorem} 
\label{thm=real}
All possible values of $\underline{d}$ and of the signature
of homogeneous $WACI$'s with formal dimension $m=4$ or $8$ may
be realized by smooth manifolds.
\end{theorem}

The rest of this section will be devoted to the proof of the 
above theorem. The case $m=4$ is easy.
 
\begin{lemma}
\label{lem:four}
Theorem~\ref{thm=real} is true for $m=4$.
\end{lemma}

\begin{proof}
In case $(I_4)$, $\mid \sigma \mid =2$ or $0$, realized by
$H^*(\C \PP^2 \# \C \PP^2 , \Q) = \Q[x_1, x_2]/(x^2_1 -x^2_2, x_1 x_2)$,
and $H^*(\C \PP^1 \times \C \PP^1 , \Q) = \Q[x_1, x_2]/(x^2_1, x^2_2)$
respectively.

In case $(II_4)$, $\mid \sigma \mid =1$, realized by
$H^*(\C \PP^2 , \Q) = \Q[x]/(x^3)$.
\end{proof}

We move now to the case $m=8$. In the next lemma, we take care of
the subcases where the desired manifold may be obtained from known
examples, by taking products and connected sums.

\begin{lemma}
\label{lem=pr8}
Theorem~\ref{thm=real} is true for $m=8$, in all cases different from
$(IV_8)$, $\sigma =3$, and $(I_8)$, $\sigma =2$ or $6$. 
\end{lemma}

\begin{proof}
If $\sigma =0$ or $1$, one may use products of complex projective spaces,
as explained before.

{\em Case} $(I_8)$, $\sigma =4$: 
$H^*((\C \PP^2 \# \C \PP^2) \times (\C \PP^2 \# \C \PP^2) , \Q)$ 
is equal to
\[ 
\Q[x_1, x_2, y_1,y_2]/(x^2_1 -x^2_2, x_1 x_2, y^2_1 -y^2_2, y_1 y_2)\, .
\]

{\em Case} $(II_8)$, $\sigma =2$: 
$H^*((\C \PP^2 \# \C \PP^2) \times \C \PP^2, \Q) 
= \Q[x_1, x_2, x_3]/(x_1^2 -x_2^2, x_1 x_2, x_3^3)$.

{\em Case} $(II_8)$, $\sigma =4$: 
$\A = \Q [x_1, x_2, x_3]/ (x_1^2 -x_3^2, x_2^2 -x_3^2, x_1 x_2 x_3)$
is a smoothable homogeneous $WACI$, with signature $4$; see
\cite[Proposition 4.6]{Pp}. 

{\em Case} $(III_8)$, $\sigma =2$: 
$H^*(\C \PP^4 \# \C \PP^4 , \Q) = \Q[x_1, x_2]/(x^4_1 -x^4_2, x_1 x_2)$.
\end{proof}

To check the remaining cases, we will use the integrality test
from Theorem~\ref{thm=test}. Case $(IV_8)$ will follow from 
the analysis of the family below.

\begin{example}
\label{ex=c1}
Let $\A (c)$, $c\in \Q$, be the graded algebra
\[
\Q [x,y]/ (f_1 =x^3 -xy^2\, ,\, f_2 =y^3 -c x^2 y)\, ,
\]
with $x$ and $y$ of degree $2$.

It is immediate to see that $\A (c)$ is a $WACI$
(homogeneous, belonging to case $(IV_8)$) precisely
when $\{ f_1 =f_2 =0 \}= \{0 \}$ (over $\C$), that is, 
if and only if $c\ne 1$. 
\end{example}

The next lemma completes the proof of Theorem~\ref{thm=real},
case $(IV_8)$, and illustrates the arithmetic behind
integrality condition \eqref{eq:z}.

\begin{lemma}
\label{lem=f1}
Let $\{ \A (c) \}_{c\ne 1}$ be the above $WACI$ family. Then:
\begin{enumerate}
\item \label{eq=a2}
The absolute value of the signature of $\A (c)$ is
$2+\epsilon$, where $\epsilon =\sgn (c-1)$.
\item \label{eq=a1}
$\A (c)$ is smoothable $\eqv \mid c-1 \mid$ is a sum 
of two rational squares. 
\end{enumerate}
\end{lemma}

\begin{proof}
Part \eqref{eq=a2}. It is readily checked that the matrix of 
the Poincar\' e quadratic form on $\A^4 (c)$, with respect to
the basis $\{ xy, x^2, x^2 -y^2 \}$ and the orientation 
$\omega = (c-1) x^4$, is $A(c)=\diag 
(\frac{1}{c-1}, \frac{1}{c-1}, 1)$.
Clearly, the signature of $\A(c)$ is as asserted.

Part \eqref{eq=a1}. To decide the smoothability of $\A(c)$,
we will use Theorem~\ref{thm=test}\eqref{eq=t1}. Obviously,
the orientation $\omega$ satisfies the required discriminant
property. Therefore, $\A(c)$ is smoothable if and only if
$(\frac{1}{c-1}, \frac{1}{c-1})_p =1$, at all odd primes.
This equivalent (\cite[III.1]{Se}) with 
$(\frac{1}{\epsilon (c-1)}, \frac{1}{\epsilon (c-1)})_p =1$, 
at all odd primes and also at $\infty$. By Hilbert's theorem
(see~\cite[III.2]{Se}), this is further
equivalent with 
$(\frac{1}{\epsilon (c-1)}, \frac{1}{\epsilon (c-1)})_p =1$, 
at all primes and also at $\infty$.

The definition of Hilbert symbols (\cite[III.1]{Se}) and
the Hasse--Minkowski theorem (\cite[IV.3]{Se}) together
imply that this happens if and only if $\epsilon (c-1)$
is a sum of two rational squares. The proof of 
Part~\eqref{eq=a1} is complete.
\end{proof}

The last case of Theorem~\ref{thm=real} 
($(I_8)$, $\sigma =2$ or $6$) will be covered by 
analysing a second family.

\begin{example}
\label{ex:par}
Let us consider the family of graded algebras 
$\{ \B (c) \}_{c\in \Q}$, with weight $2$ generators,
$\{ x_i \}_{1\le i \le 4}$, and defining relations
\begin{equation}
\label{eq:c}
\begin{cases}
x^2_i -x^2_4 \, ,\quad \text{for} \quad i\le 3\, ,\\
\sum_{1\le i<j\le 4} x_i x_j -c x^2_4\, .
\end{cases}
\end{equation}

It is easy to see that \eqref{eq:c} defines a $WACI$
(homogeneous, belonging to case $(I_8)$) 
if and only if $c\ne -2, 0, 6$. For $c=-1$, 
\eqref{eq:c} defines the signature $6$ algebra from
\cite[p.24]{EL} (which is not smoothable, by 
Lemma~\ref{lem=f2} below).
\end{example}

The next lemma completes the proof of Theorem~\ref{thm=real}.
For its proof, we will resort to the integrality test from
Theorem~\ref{thm=test}\eqref{eq=t2}. At this point, it seems
worthwhile pointing out that Part~\eqref{eq=b3} of the lemma
provides an interesting family of examples, where property
\eqref{eq:z} may be decided using only the (simple)
discriminant obstruction. This simple behaviour cannot be
expected, in general; see Remark~\ref{rk=w}. 

\begin{lemma}
\label{lem=f2}
Let $\{ \B (c) \}_{c\ne -2,0,6}$ be the above $WACI$ family. Then:
\begin{enumerate}
\item \label{eq=b1}
The signature of $\B (c)$ is $0, \pm 2$ or $\pm 6$.
\item \label{eq=b3}
$\B (c)$ is smoothable if and only if 
$\mid (c-6)(c+2) \mid \in \Q^{*2}$.
\item \label{eq=b2}
For $c=-3, 2$ and $-\frac{2}{5}$, the algebra
$\B (c)$ is smoothable, with signature
$0, 2$ and $6$ respectively.
\end{enumerate}
\end{lemma}

\begin{proof}
Part~\eqref{eq=b1}. Our first task is to find a $\Q$--basis of
$\B^4(c)$, and an orientation $\omega$, with respect to which
the matrix of the Poincar\' e inner product is diagonal. We will
begin with the basis $\{ x_i x_j \}_{1\le i<j\le 4}$. Set
$y:= x_i^2\in \B^4(c)$, $1\le i\le 4$. We claim that
\begin{equation}
\label{eq=m1}
x_1 x_2 x_3 x_4=\frac{c^2 -4c -6}{6} y^2\, ,
\end{equation}
and
\begin{equation}
\label{eq=m2}
y x_i x_j= \frac{c}{6} y^2\, ,\quad \text{for} \quad
1\le i<j\le 4\, ,
\end{equation}
so we may take $\omega =y^2$. Indeed, we infer from 
\eqref{eq:c} that 
\begin{equation}
\label{eq=m3}
cy x_1 x_2 =cy x_3 x_4=y^2 +y(x_1 +x_2)(x_3 +x_4)+x_1 x_2 x_3 x_4 \, .
\end{equation}
Adding all relations of type \eqref{eq=m3}, we get \eqref{eq=m1}.
Using \eqref{eq=m1}, \eqref{eq=m2} follows from \eqref{eq=m3}.

Consider now the following basis in the middle dimension $4$:
$e_1=x_1x_2-x_3x_4, e_2=x_1x_4-x_3x_2, e_3=x_1x_3-x_2x_4, e_4=x_1x_2+x_3x_4, 
e_5=x_1x_4+x_3x_2, e_6=x_1x_3+x_2x_4$, and rescale the orientation to
$\frac{(6-c)(c+2)}{3}y^2$. 

With respect to these data, the intersection form is given by the 
following matrix:
\[
B_1(c)=
\begin{pmatrix}

1&0&0&0&0&0\\
0&1&0&0&0&0\\
0&0&1&0&0&0\\
0&0&0&a&b&b\\
0&0&0&b&a&b\\
0&0&0&b&b&a
\end{pmatrix}
\]
where $a=\frac{c(c-4)}{(6-c)(c+2)}, b=\frac{2c}{(6-c)(c+2)}$; use
\eqref{eq=m1} and \eqref{eq=m2}.  

Firstly note that if $c=4$ i.e. $a=0$, the determinant of $B_1$ 
is positive, so clearly this case cannot produce $\sigma =4$.

Let us now consider the case $c\neq 4$, i.e. $a\neq 0$.

Considering a new basis: 
$f_i=e_i, i \leq 4, f_5=e_5-e_6, f_6= -2be_4+ae_5+ae_6$, 
we obtain a new matrix:
\[
B_2(c)=
\begin{pmatrix}

1&0&0&0&0&0\\
0&1&0&0&0&0\\
0&0&1&0&0&0\\
0&0&0&a&0&0\\
0&0&0&0&2(a-b)&0\\
0&0&0&0&0&2a(a-b)(a+2b)
\end{pmatrix}
\]

Rescaling our last basis to: 
$ g_i=f_i, i\leq 5, g_6=\frac{(6-c)(c+2)}{c^2}f_6$, 
we finally obtain the matrix:
\begin{equation}
\label{eq=diag}
B_3(c)=
\begin{pmatrix}

1&0&0&0&0&0\\
0&1&0&0&0&0\\
0&0&1&0&0&0\\
0&0&0&\frac{c(c-4)}{(6-c)(c+2)}&0&0\\
0&0&0&0&\frac{-2c}{c+2}&0\\
0&0&0&0&0&\frac{-2(c-4)}{c+2}
\end{pmatrix}
\end{equation}

Our assertion on signature from Part~\eqref{eq=b1}
easily follows by examining the distribution of signs
on the diagonal of the matrix $B_3(c)$. 

Part~\eqref{eq=b2}. Follows from Part~\eqref{eq=b3}.

Part~\eqref{eq=b3}. The algebra $\B(c)$ is smoothable 
if and only if it verifies the integrality test from 
Theorem~\ref{thm=test}~\eqref{eq=t2}. The discriminant 
may be computed as $\det B_1(c) \equiv (6-c)(c+2)$
(modulo $\Q^{*2}$). This shows that we may assume from now on 
$c\ne 4$, and use the matrix $B_3(c)$. Set 
$\epsilon :=\sgn ((6-c)(c+2))$, and note that $\epsilon =+1$
(respectively $\epsilon =-1$) corresponds to the case
$(\bullet \bullet)$ (respectively $(\bullet)$). We have to
show that, in both cases, the property
$\mid (6-c)(c+2) \mid \in \Q^{*2}$ implies the restrictions on
the local invariants of $B:= B_3(c)$ from Theorem~\ref{thm=test}. Set
\begin{equation} 
\label{eq=3l}
\lambda_1=\frac{-2c}{c+2},\quad \lambda_2=\frac{-2(c-4)}{c+2},
\quad \text{and} \quad \lambda_3=\frac{c(c-4)}{(6-c)(c+2)}\, ,
\end{equation}
and note that $\lambda_3 \equiv \epsilon \lambda_1 \lambda_2$
(modulo $\Q^{*2}$), by our assumption on the discriminant. By 
elementary manipulations with Hilbert symbols, we infer that
\begin{equation}
\label{eq=ecases}
\varepsilon_p(B)=
\begin{cases}
(\lambda_1 ,\lambda_2)_p\, , & \text{if} \quad \epsilon =-1\, ;\\
(\lambda_1 ,\lambda_2)_p\cdot (\lambda_1 ,\lambda_1)_p \cdot
(\lambda_2 ,\lambda_2)_p\, ,  & \text{if} \quad \epsilon =+1\, .
\end{cases}
\end{equation}

The discriminant condition means that
\[
\epsilon (6-c) =\frac{t^2}{s^2} (c+2)\, ,
\]
where $t$ and $s$ are relatively prime integers. Solve for $c$
and substitute in \eqref{eq=3l}, to obtain the following values
(modulo $\Q^{*2}$) for $\lambda_{1,2}$:
\begin{equation}
\label{eq=st}
\begin{cases}
\lambda_1= 2\epsilon t^2- 6s^2\, ,\\
\lambda_2= 6\epsilon t^2- 2s^2\, .
\end{cases}
\end{equation}

We are going to compute the Hilbert symbols appearing 
in \eqref{eq=ecases}, in terms of Legendre symbols; see
\cite[I.3 and Theorem III.1]{Se}. To do this, write
\begin{equation}
\label{eq=ab}
\begin{cases}
\lambda_1=p^{\alpha}u \, ,\\
\lambda_2=p^{\beta}v \, ,
\end{cases}
\end{equation}
where $\alpha, \beta \in \N$, $u,v\in \Z$, and
$u,v \not\equiv 0 (p)$.

To finish our proof, we are going to show that
$\varepsilon_p(B)=1$, $\forall p\equiv 1 (2)$
(when $\epsilon=-1$), and
$\varepsilon_p(B)=1$, $\forall p\equiv 1 (4)$
(when $\epsilon=+1$).

Several cases may appear in \eqref{eq=ab}. If
$\alpha =\beta =0$, then plainly 
$(\lambda_1 ,\lambda_2)_p= (\lambda_1 ,\lambda_1)_p=
(\lambda_2 ,\lambda_2)_p=1$, at all odd primes $p$.
The case $\alpha, \beta>0$ cannot occur, since this 
would imply (see \eqref{eq=st}) that $s\equiv t\equiv 0 (p)$.
The remaining cases ($\alpha =0, \beta>0$ and $\alpha>0, \beta=0$)
may be settled as follows.

For $\alpha =0, \beta>0$, one knows (\cite[Theorem III.1]{Se})
that $(\lambda_1, \lambda_2)_p=\Big ( \frac{\lambda_1}{p} \Big )^{\beta}$.
Since $2s^2 \equiv 6\epsilon t^2\, (p)$, 
$\lambda_1\equiv -16\epsilon t^2\, (p)$. Therefore,
$(\lambda_1, \lambda_2)_p= (-\epsilon)^{\beta \varepsilon (p)}$.
Similarly, for $\alpha>0, \beta=0$, one has
$(\lambda_1, \lambda_2)_p=\Big ( \frac{\lambda_2}{p} \Big )^{\alpha}$,
with $\lambda_2 \equiv 16 s^2\, (p)$, hence 
$(\lambda_1, \lambda_2)_p=1$. By \eqref{eq=ecases}, this 
completes our proof, when $\epsilon =-1$.

Assume now $\epsilon =+1$. In this last case, we will
also need $(\lambda_1, \lambda_1)_p= (-1)^{\alpha \varepsilon (p)}$,
and $(\lambda_2, \lambda_2)_p= (-1)^{\beta \varepsilon (p)}$.
When $p\equiv 1(4)$, both $(\lambda_1, \lambda_1)_p$ and
$(\lambda_2, \lambda_2)_p$ are $1$, which completes our proof 
(see \eqref{eq=ecases}).

\end{proof}


\end{document}